\documentclass[a4paper,11pt]{article}
\usepackage{amsmath}
\usepackage{amsfonts}
\usepackage{amsthm}
\usepackage{amssymb}
\usepackage{latexsym}
\usepackage{fancyhdr}
\usepackage[latin1]{inputenc}
\usepackage{makeidx}

\newcommand{\bbt}{\mathbb{T}}
\newcommand{\bbd}{\mathbb{D}}

\newcommand{\bbz}{\mathbb{Z}}
\newcommand{\bbr}{\mathbb{R}}
\newcommand{\bbc}{\mathbb{C}}
\newcommand{\calc}{\mathcal{C}}
\newcommand{\calh}{\mathcal{H}}
\newcommand{\dsp}{\displaystyle}
\newcommand{\e}{\varepsilon}

\newtheorem{defi}{D\'efinition}[section]

\newtheorem{theo}{Th\'eor\`eme}[section]
\newtheorem{prop}[theo]{Proposition}

\newtheorem{lemme}[theo]{Lemme}

\makeindex

\begin{document}

\title{Id\'eaux ferm\'es de certaines alg\`ebres de Beurling \\
et application aux op\'erateurs \`a spectre d\'enombrable}
\date{}
\maketitle

\begin{abstract}
We denote by $\bbt$ the unit circle and by $\bbd$ the unit disc of $\bbc$. Let $s$ be a non-negative real and $\omega$ a weight such that $\omega(n) = (1+n)^{s} \quad (n \geq 0)$ and such that the sequence $\dsp \Big( \frac{\omega(-n)}{(1+n)^{s}} \Big)_{n \geq 0}$ is non-decreasing. We define the Banach algebra
$$
A_{\omega}(\bbt) = \Big\{ f \in \calc(\bbt) : \, \big\| f \big\|_{\omega} = \sum_{n = -\infty}^{+\infty} | \hat{f}(n) | \omega(n) < +\infty \Big\},
$$
If $I$ is a closed ideal of $A_{\omega}(\bbt)$, we set $h^{0}(I) = \Big\{ z \in \bbt : \, f(z) = 0 \quad (f \in I) \Big\}$. We describe here all closed ideals $I$ of $A_{\omega}(\bbt)$ such that $h^{0}(I)$ is at most countable. A similar result is obtained for closed ideals of the algebra $A_{s}^{+}(\bbt) = \Big\{ f \in A_{\omega}(\bbt) : \, \widehat{f}(n) = 0 \quad (n<0) \Big\}$ without inner factor. \\
Then, we use this description to establish a link between operators with countable spectrum and interpolating sets for $\textrm{{\LARGE $a$}}^{\infty}$, the space of infinitely differentiable functions in the closed unit disc $\overline{\bbd}$ and holomorphic in $\bbd$. 
\end{abstract}

   \section{Introduction}

\thispagestyle{fancy}
\lfoot{AMS Subject Classification 2000: 46J20, 47A30, 30H05 } 
\cfoot{\ \\1}
\renewcommand{\headrulewidth}{0pt}
\renewcommand{\footrulewidth}{0pt}

On note $\bbt$ le cercle unit\'e et $\bbd$ le disque unit\'e de $\bbc$. Pour un entier $p \geq 0$, on note $\calc^{p}(\bbt)$ l'alg\`ebre des fonctions $p$ fois contin\^ument d\'erivables sur $\bbt$. Si $\omega = \big( \omega(n) \big)_{n \in \bbz}$ est une suite de r\'eels strictement positifs, on d\'efinit $A_{\omega}(\bbt)$ par
$$
A_{\omega}(\bbt) = \Big\{ f \in \calc^{0}(\bbt) : \, \big\| f \big\|_{\omega} = \sum_{n = -\infty}^{+\infty} | \hat{f}(n) | \omega(n) < +\infty \Big\},
$$
o\`u $\hat{f}(n)$ d\'esigne le $\textrm{n}^{\textrm{i\`eme}}$ coefficient de Fourier de $f$. On dit qu'une suite $\omega = \big( \omega(n) \big)_{n \in \bbz}$ est un poids si $\omega(n) \geq 1$ et $\omega(m+n) \leq \omega(m) \omega(n)$ pour tout $m,n$ dans $\bbz$. Si $\omega$ est un poids, $A_{\omega}(\bbt)$ muni de la norme $\| \, \|_{\omega}$ est une alg\`ebre de Banach. Elle est r\'eguli\`ere si et seulement si $\dsp \sum_{n=-\infty}^{+\infty} \frac{\log \omega(n)}{1+n^{2}} < +\infty$ (voir \cite{Katz}, ex.7 p.118). Soit $p$ un entier positif tel que  $A_{\omega}(\bbt) \subset \calc^{p}(\bbt)$. Si $I$ est un id\'eal ferm\'e de $A_{\omega}(\bbt)$, on pose
$$
h^{k}(I) = \Big\{ z \in \bbt : \, f(z) = \ldots = f^{(k)}(z) = 0 \quad (f \in I) \Big\} \quad \big( k \in \{0, \ldots, p\} \big).
$$

Soit $s$ un r\'eel positif, on note $[s]$ sa partie enti\`ere. On dira qu'une suite de r\'eels strictement positifs $\omega = \big( \omega(n) \big)_{n \in \bbz}$ v\'erifie la condition $(W_{s})$ si
$$ 
\left\{ \begin{array}{lll}
\omega(n) = (1+n)^{s} \qquad (n \geq 0) \\   
\textrm{la suite } \, \dsp \Big( \frac{\omega(-n)}{(1+n)^{s}} \Big)_{n \geq 0} \textrm{ est croissante},
\end{array} \right. \quad \eqno (W_{s}) 
$$
et qu'il v\'erifie la condition $(A)$ si
$$ \label{condA}
\qquad \omega(-n) = O \big( e^{\e \sqrt{n}} \big) \,\, (n \rightarrow +\infty), \, \textrm{ pour tout } \e > 0. \eqno (A)
$$ 
Soit $\omega = \big( \omega(n) \big)_{n \in \bbz}$ un poids v\'erifiant les conditions $(W_{s})$ et $(A)$. On a alors l'inclusion $A_{\omega}(\bbt) \subset \calc^{[s]}(\bbt)$, et $A_{\omega}(\bbt)$ est une alg\`ebre de Banach r\'eguli\`ere. Dans le cas particulier o\`u le poids $\omega$ est d\'efini par $\omega(n) = (1+|n|)^{s}$ pour tout $n \in \bbz$, avec $s$ un r\'eel positif, on notera $\big( A_{s}(\bbt), \| \, \|_{s} \big)$ l'alg\`ebre $\big( A_{\omega}(\bbt), \| \, \|_{\omega} \big)$. On remarque que $A_{0}(\bbt) = A(\bbt)$ n'est rien d'autre que l'alg\`ebre de Wiener. On posera \'egalement
$$
A_{s}^{+}(\bbt) = \Big\{ f \in A_{s}(\bbt) : \, \widehat{f}(n) = 0 \quad (n<0) \Big\}.
$$
Pour $f \in A_{s}^{+}(\bbt)$, on note $S(f)$ son facteur int\'erieur et on pose
$$
Z_{+}^{k}(f) = \Big\{ z \in \overline{\bbd} : \, f(z) = \ldots = f^{(k)}(z) = 0 \Big\} \quad \big( k \in \{0, \ldots, [s]\} \big).
$$
Si $I$ est un id\'eal ferm\'e de $A_{s}^{+}(\bbt)$, on note $S_{I}$ son facteur int\'erieur, c'est-\`a-dire le plus grand diviseur int\'erieur commun \`a tous les \'el\'ements de $I$ non nuls (voir \cite{Hoff} p.85), et on pose $\dsp h_{+}^{k}(I) = \bigcap_{f \in I} Z_{+}^{k}(f)$.

Nous d\'ecrivons dans un premier temps tous les id\'eaux ferm\'es $I$ de $A_{\omega}(\bbt)$ lorsque $h^{0}(I)$ est au plus d\'enombrable et lorsque le poids $\omega$ v\'erifie les conditions $(W_{s})$ et $(A)$. Plus pr\'ecis\'ement, nous montrons (th\'eor\`eme \ref{ideauxdeaw}) que, sous ces conditions, on a
\begin{eqnarray*} 
I = \Big\{ f \in A_{\omega}(\bbt) : \, f^{(j)}(z)=0 \textrm{ sur } h^{j}(I) \quad (0 \leq j \leq [s]) \Big\}. 
\end{eqnarray*}
Dans le cas $s=0$, nous retrouvons ainsi un r\'esultat connu (voir \cite{Zar1}). Nous d\'eduisons alors de ce r\'esultat une caract\'erisation des id\'eaux ferm\'es de $A_{s}^{+}(\bbt)$ sans facteur int\'erieur (c'est-\`a-dire tel que $S_{I}=1$) tels que $h_{+}^{0}(I)$ est au plus d\'enombrable. On notera $\calh^{\infty}(\bbd)$ l'alg\`ebre des fonctions holomorphes et born\'ees dans $\bbd$ et si $E_{[s]} \subset \ldots \subset E_{0}$ sont des ferm\'es de $\bbt$ et $S$ une fonction int\'erieure, on d\'efinit 
$$
I \big( S ; \, E_{0}, \ldots, E_{[s]} \big) = \Big\{ f \in A_{s}^{+}(\bbt) : \, S \big| S(f), \, E_{0} \subset Z_{+}^{0}(f) \cap \bbt, \ldots, E_{[s]} \subset Z_{+}^{[s]}(f) \cap \bbt \Big\},
$$
o\`u $S \big| S(f)$  signifie que $S$ divise $S(f)$, c'est-\`a-dire que $\dsp \frac{S(f)}{S}$ est dans $\calh^{\infty}(\bbd)$. Nous montrons que si $I$ est un id\'eal ferm\'e de $A_{s}^{+}(\bbt)$ sans facteur int\'erieur tel que $h_{+}^{0}(I)$ est au plus d\'enombrable, alors
$$
I = I \big( 1 ; \, h_{+}^{0}(I) \cap \bbt, \ldots, h_{+}^{[s]}(I) \cap \bbt \big).
$$
Notons que dans le cas $s=0$, les id\'eaux ferm\'es $I$ de $A^{+}(\bbt) = A_{0}^{+}(\bbt)$ ont \'et\'e caract\'eris\'es par J. P. Kahane dans \cite{Kaha2} lorsque $h_{+}^{0}(I)$ est fini, par C. Bennett et J. E. Gilbert dans \cite{BeGi} lorsque $h_{+}^{0}(I)$ est d\'enombrable, et enfin par J. Esterle, E. Strouse et F. Zouakia dans \cite{EStZo2} lorsque $h_{+}^{0}(I)$ est le Cantor triadique. D'autre part, J. Esterle a construit un id\'eal ferm\'e $I$ de $A^{+}(\bbt)$ tel que $I \neq I \big( S_{I} ; h_{+}^{0}(I) \big)$, d\'emontrant ainsi que la conjecture de C. Bennett et J. E. Gilbert dans \cite{BeGi} est fausse. \\

Dans un deuxi\`eme temps, nous allons utiliser ce r\'esultat pour \'etudier le comportement de certains op\'erateurs \`a spectre d\'enombrable et inclus dans $\bbt$. Soit $E$ un ferm\'e de $\bbt$ et $s$, $t$ deux r\'eels positifs ou nuls. On d\'esignera par $P(s,t,E)$ la propri\'et\'e suivante: tout op\'erateur $T$ inversible sur un espace de Banach tel que $\textrm{Sp} \, T \subset E$ et qui v\'erifie les conditions:
\begin{eqnarray} 
& & \big\| T^{n} \big\| = O(n^{s}) \,\, (n \rightarrow +\infty) \label{introCp} \\
& & \big\| T^{-n} \big\| = O(e^{\e \sqrt{n}}) \,\, (n \rightarrow +\infty), \, \textrm{ pour tout } \e > 0, \label{introA} 
\end{eqnarray}
v\'erifie \'egalement la propri\'et\'e plus forte 
\begin{eqnarray}
\big\| T^{-n} \big\| = O(n^{t}) \, (n \rightarrow +\infty). \qquad\qquad\quad\, \label{introDq}
\end{eqnarray}
M. Zarrabi a montr\'e dans \cite{Zar1} (th\'eor\`eme 3.1 et remarque 2.a) qu'un ferm\'e $E$ de $\bbt$ v\'erifie la propri\'et\'e $P(0,0,E)$ si et seulement si $E$ est d\'enombrable. \\

Nous nous proposons d'\'etudier la propri\'et\'e $P(s,t,E)$ pour n'importe quel r\'eel $s \geq 0$. On dira qu'un ferm\'e $E$ de $\bbt$ v\'erifie la condition de Carleson si 
$$
\int_{0}^{2\pi} \log^{+} \frac{1}{d(e^{it},E)} \mathrm{d} t < +\infty, \eqno (C)
$$
et qu'il v\'erifie la condition (ATW) s'il existe des constantes $C_{1},C_{2} > 0$ telles que
$$
\frac{1}{|L|} \int_{L} \log^{+} \frac{1}{d(e^{it},E)} \mathrm{d} t \leq C_{1}\log \frac{1}{|L|} + C_{2}\, , \quad \textrm{pour tout arc $L$ de $\bbt$}, \eqno (ATW)
$$ 
o\`u $|L|$ d\'esigne la longueur de l'arc $L$, et $d(e^{it},E)$ la distance de $e^{it}$ \`a $E$. \\
La condition (ATW) vient de \cite{ATaWi2}, o\`u les auteurs montrent que les ensembles v\'erifiant (ATW) sont les ensembles d'interpolation pour $\textrm{{\LARGE $a$}}^{\infty}$, l'espace des fonctions holomorphes dans le disque unit\'e et de classe $\calc^{\infty}$ dans $\overline{\bbd}$. \\

Nous montrons (th\'eor\`eme \ref{thetheo}) que si $E$ est un ferm\'e d\'enombrable de $\bbt$, alors les deux assertions suivantes sont \'equivalentes: \\
\hspace*{3mm} (i) $E$ v\'erifie la condition (ATW). \\
\hspace*{2mm} (ii) $E$ v\'erifie la condition (C) et pour tout r\'eel $s \geq 0$, il existe un r\'eel $t$ tel que la propri\'et\'e $P(s,t,E)$ soit v\'erifi\'ee. \\

Nous montrons ensuite que si $E$ n'est pas d\'enombrable, alors la propri\'et\'e $P(0,t,E)$ n'est v\'erifi\'ee pour aucun r\'eel $t \geq 0$. \\

   \section{Unit\'e approch\'ee pour certains id\'eaux de $A_{\omega}(\bbt)$}

Nous allons \'etablir ici que lorsque le poids $\omega$ satisfait la condition $(W_{s})$, alors l'alg\`ebre $A_{\omega}(\bbt)$ v\'erifie la condition de Ditkin analytique forte (au sens de \cite{BeGi}, p. 4). Pour cela, nous avons besoin des trois lemmes \'el\'ementaires suivants: \\

\begin{lemme} \label{inegal}
Soient $\beta$ un r\'eel tel que $0 \leq \beta < 1$ et $j$ un entier positif. Alors pour tout r\'eel $x$ tel que $0 \leq x < 1$, nous avons l'in\'egalit\'e suivante:
$$
\sum_{k=j}^{+\infty} (1+k)^{\beta} x^{k} \leq \frac{(j+1)^{\beta} x^{j}}{1-x} + \frac{x^{j+1}}{(1-x)^{\beta+1}}.
$$
\end{lemme}

\begin{proof}
Soit $x$ un r\'eel tel que $0 \leq x < 1$, posons
$$
f(x) = \sum_{k=j}^{+\infty} (1+k)^{\beta} x^{k} - \frac{(j+1)^{\beta} x^{j}}{1-x} - \frac{x^{j+1}}{(1-x)^{\beta+1}}.
$$
En d\'eveloppant en s\'erie les fonctions $\dsp x \mapsto \frac{1}{1-x}$ et $\dsp x \mapsto \frac{1}{(1-x)^{\beta+1}}$, on montre que
$$
f(x) = \sum_{k=1}^{+\infty} a_{k} x^{k+j},
$$
avec $a_{1} = (j+2)^{\beta} - (j+1)^{\beta} - 1$ et pour $k \geq 2$,
$$
a_{k} = (k+j+1)^{\beta} - (j+1)^{\beta} - \frac{(\beta+1) \ldots (\beta+k-1)}{(k-1)!}.
$$
Il est facile de voir que $a_{1} \leq 0$. Pour d\'emontrer ce lemme, il suffit donc de montrer que pour tout $k \geq 2$, $a_{k} \leq 0$. On commence par remarquer que, pour tout $k \geq 2$, la fonction $\dsp t \mapsto (k+t+1)^{\beta} - (t+1)^{\beta} - \frac{(\beta+1) \ldots (\beta+k-1)}{(k-1)!}$ est d\'ecroissante sur $[-1,+\infty)$. Par cons\'equent, pour tout $k \geq 2$, on a
\begin{eqnarray} \label{1_inegal1}
a_{k} \leq k^{\beta} - \frac{(\beta+1) \ldots (\beta+k-1)}{(k-1)!}.
\end{eqnarray}
Puisque $\beta < 1$, on a $(1+t)^{\beta} \leq 1 + \beta t$, pour tout r\'eel $t$ positif. En utilisant cette in\'egalit\'e, on obtient 
\begin{eqnarray*}
\frac{(\beta+1) \ldots (\beta+k-1)}{(k-1)!} & = & (1+\beta) \Big( 1 + \frac{\beta}{2} \Big) \ldots \Big( 1 + \frac{\beta}{k-1} \Big) \\
& \geq & (1+1)^{\beta} \Big( 1 + \frac{1}{2} \Big)^{\beta} \ldots \Big( 1 + \frac{1}{k-1} \Big)^{\beta} = k^{\beta}.
\end{eqnarray*}
On d\'eduit alors de cette in\'egalit\'e et de (\ref{1_inegal1}) que pour tout $k \geq 2$, $a_{k} \leq 0$, ce qu'il s'agissait de d\'emontrer.
\end{proof}

\begin{lemme} \label{derivee}
Soit $\omega = \big( \omega(n) \big)_{n \in \bbz}$ une suite de r\'eels strictement positifs telle que : \\
(i) les suites $\dsp \Big( \frac{\omega(-n)}{1+n} \Big)_{n \geq 0}$ et $\dsp \Big( \frac{\omega(n)}{1+n} \Big)_{n \geq 0}$ soient croissantes. \\
(ii) $\dsp 0 < \inf_{n \in \bbz} \frac{\omega(n+1)}{\omega(n)} < \sup_{n \in \bbz} \frac{\omega(n+1)}{\omega(n)} < +\infty$. \\
On note $\dsp A = \inf_{n \in \bbz} \frac{\omega(n+1)}{\omega(n)}$ et $\dsp B = \sup_{n \in \bbz} \frac{\omega(n+1)}{\omega(n)}$, et on d\'efinit la suite $\omega_{1} = \big( \omega_{1}(n) \big)_{n \in \bbz}$ par 
$$
\omega_{1}(n) = \frac{\omega(n)}{1+|n|} \qquad (n \in \bbz).
$$
Alors on a les deux propri\'et\'es suivantes: \\
1) On a la double in\'egalit\'e $A \big\| f' \big\|_{\omega_{1}} \leq \big\| f \big\|_{\omega} \leq \big| \hat{f}(0) \big| \, \omega(0) + 3B \big\| f' \big\|_{\omega_{1}}$. \\
2) Soit $f$ dans $A_{\omega}(\bbt)$ qui s'annule au point $1$, alors $\dsp \frac{f}{\alpha-1}$ est dans $A_{\omega_{1}}(\bbt)$, o\`u $\alpha: z \mapsto z$.
\end{lemme}

\begin{proof}
1) Cela d\'ecoule imm\'ediatement de la relation $\widehat{f'}(n) = (n+1) \hat{f}(n+1) \,\, (n \in \bbz)$. \\
2) Soit $f$ dans $A_{\omega}(\bbt)$ qui s'annule au point $1$, on \'ecrit $f = \sum \limits_{n = -\infty}^{+\infty}  \hat{f}(n) (\alpha^{n} - 1)$, d'o\`u
$$
\frac{f}{\alpha-1} = - \sum_{n = -\infty}^{-1} \hat{f}(n) \big( \alpha^{n} + \ldots + \alpha^{-1}\big) +  \sum_{n = 1}^{+\infty} \hat{f}(n) \big( 1 + \ldots + \alpha^{n-1} \big). 
$$
Ainsi
\begin{eqnarray} \label{lemme21}
\Big\| \frac{f}{\alpha-1} \Big\|_{\omega_{1}} \, = \, \sum_{n = -\infty}^{-1} \big| \hat{f}(n) \big| \Big( \omega_{1}(n) + \ldots + \omega_{1}(-1) \Big)\, + \, \sum_{n = 1}^{+\infty} \big| \hat{f}(n) \big| \Big( \omega_{1}(0) + \ldots + \omega_{1}(n-1) \Big). \quad 
\end{eqnarray}
Comme les suites $\big( \omega_{1}(-n) \big)_{n \geq 0}$ et $\big( \omega_{1}(n) \big)_{n \geq 0}$ sont croissantes, on d\'eduit de (\ref{lemme21}) que
\begin{eqnarray*}
\Big\| \frac{f}{\alpha-1} \Big\|_{\omega_{1}} & \leq & \sum_{n = -\infty}^{-1} \big| \hat{f}(n) \big| \, |n| \, \omega_{1}(n) \, + \, \sum_{n = 1}^{+\infty} \big| \hat{f}(n) \big| \, n \, \omega_{1}(n) \\
& \leq & \sum_{n = -\infty}^{+\infty} \big| \hat{f}(n) \big| \, \omega(n) < +\infty.
\end{eqnarray*}
\end{proof}

\begin{lemme} \label{unite}
Soient $0 \leq \beta < 1$ un r\'eel et $\omega = \big( \omega(n) \big)_{n \in \bbz}$ une suite v\'erifiant la condition $(W_{\beta})$. On d\'efinit une suite de fonctions de $A_{\omega}(\bbt)$ par  
\begin{eqnarray}  \label{unite0}
e_{n} = \frac{\alpha-1}{\alpha-1-\frac{1}{n}} \qquad (n \geq 1),
\end{eqnarray}
o\`u $\alpha$ est la fonction $z \mapsto z$. Alors pour toute fonction $f$ dans $A_{\omega}(\bbt)$ telle que $f(1)=0$, on a 
\begin{eqnarray} 
& & \big\| (e_{n}-1)f \big\|_{\omega} \leq 3 \big\| f \big\|_{\omega} \label{unite1} \\
& & \lim_{n \rightarrow +\infty} \big\| (e_{n}-1)f \big\|_{\omega}  = 0. \label{unite2}
\end{eqnarray}
\end{lemme}

\begin{proof}
Soit $f$ dans $A_{\omega}(\bbt)$ telle que $f(1)=0$. En \'ecrivant $f = f - f(1) = \sum \limits_{j = -\infty}^{+\infty} \hat{f}(j) (\alpha^{j} - 1)$, on voit que (\ref{unite1}) sera d\'emontr\'e si on v\'erifie que 
$$
\big\| (e_{n}-1) (\alpha^{j}-1) \big\|_{\omega} \leq 3 \, \omega(j) \quad (j \in \bbz).
$$
Soit $n \geq 1$. Un simple calcul montre que
$$
e_{n} = 1 - \frac{1}{n+1} \sum_{k = 0}^{+\infty} \Big( \frac{n}{n+1} \Big)^{k} \alpha^{k},
$$
et
$$
(e_{n}-1) (\alpha^{j}-1) = - \frac{1}{n+1} \sum_{k = 0}^{+\infty} \Big( \frac{n}{n+1} \Big)^{k} \alpha^{j+k} + \frac{1}{n+1} \sum_{k = 0}^{+\infty} \Big( \frac{n}{n+1} \Big)^{k} \alpha^{k}.
$$
Supposons que $j \geq 0$. On a
$$
\big\| (e_{n}-1) (\alpha^{j}-1) \big\|_{\omega} = \frac{1}{n+1} \sum_{k=0}^{j-1} \Big( \frac{n}{n+1} \Big)^{k} (1+k)^{\beta} + \frac{( \frac{n+1}{n})^{j}-1}{n+1} \sum_{k=j}^{+\infty} \Big( \frac{n}{n+1} \Big)^{k} (1+k)^{\beta}.
$$
On a alors l'in\'egalit\'e suivante:
\begin{eqnarray} \label{unite001}
\Big( \frac{n+1}{n} \Big)^{j} - 1 \leq \min{\big( 1,\frac{j}{n} \big)} \Big( \frac{n+1}{n} \Big)^{j}.
\end{eqnarray}
En effet, $\dsp \Big( \frac{n+1}{n} \Big)^{j} - 1 \leq \Big( \frac{n+1}{n} \Big)^{j}$, ce qui prouve l'in\'egalit\'e ci-dessus dans le cas o\`u $j \geq n$. Dans le cas o\`u $j \leq n$ et $j \neq 0$, on d\'eduit de l'in\'egalit\'e des accroissements finis que
$$
\Big( \frac{n+1}{n} \Big)^{j} - 1 \, \leq \, \frac{j}{n} \Big( \frac{n+1}{n} \Big)^{j-1} \, = \, \frac{j}{n+1} \Big( \frac{n+1}{n} \Big)^{j}.
$$
Par cons\'equent, dans le cas $j \leq n$, on a $ \dsp
\Big( \frac{n+1}{n} \Big)^{j} - 1 \leq \frac{j}{n} \Big( \frac{n+1}{n} \Big)^{j}$ (le cas $j=0$ \'etant \'evident). On d\'eduit alors de l'in\'egalit\'e (\ref{unite001}) et du lemme \ref{inegal} appliqu\'e pour $\dsp x=\frac{n}{n+1}$, la majoration suivante:
\begin{eqnarray} \label{unite01}
\big\| (e_{n}-1) (\alpha^{j}-1) \big\|_{\omega} \leq \frac{1}{n+1} \sum_{k=0}^{j-1} \Big( \frac{n}{n+1} \Big)^{k} (1+k)^{\beta} + \min{\big( 1,\frac{j}{n} \big)} \Big( (j+1)^{\beta} + (n+1)^{\beta} \Big). 
\end{eqnarray}
En observant maintenant que $\dsp \frac{1}{n+1} \sum_{k=0}^{j-1} \Big( \frac{n}{n+1} \Big)^{k} (1+k)^{\beta} \leq (1+j)^{\beta}$, et en distinguant les cas $j \leq n$ et $j \geq n+1$, on obtient
$$
\big\| (e_{n}-1) (\alpha^{j}-1) \big\|_{\omega} \leq 3 \, (1+j)^{\beta} = 3 \, \omega(j).
$$
Supposons maintenant que $j \leq -1$. On a
$$
\big\| (e_{n}-1) (\alpha^{j}-1) \big\|_{\omega} = \frac{1}{n+1} \sum_{k=j}^{-1} \Big( \frac{n}{n+1} \Big)^{k-j} \omega(k) + \frac{1-(\frac{n}{n+1})^{-j}}{n+1} \sum_{k=0}^{+\infty} \Big( \frac{n}{n+1} \Big)^{k} (1+k)^{\beta}.
$$
Avec les m\^emes arguments que ceux utilis\'es pr\'ec\'edemment, on obtient
\begin{eqnarray} \label{unite02}
\big\| (e_{n}-1) (\alpha^{j}-1) \big\|_{\omega} & \leq & \frac{1}{n+1} \sum_{k=j}^{-1} \Big( \frac{n}{n+1} \Big)^{k-j} \omega(k) + 2 \min{\big( 1, -\frac{j}{n} \big)} (n+1)^{\beta}
\end{eqnarray}
Gr\^ace aux hypoth\`eses faites sur le poids $\omega$, on a $(1+|j|)^{\beta} \leq \omega(j)$ et $\omega(k) \leq \omega(j)$ si $j \leq k \leq -1$. On en d\'eduit alors que
\begin{eqnarray*}
\big\| (e_{n}-1) (\alpha^{j}-1) \big\|_{\omega} & \leq & \omega(j) + 2 (1+|j|)^{\beta} \\
& \leq & 3 \, \omega(j).
\end{eqnarray*}
On vient donc de d\'emontrer que pour tout $j \in \bbz$,
$$
\big\| (e_{n}-1) (\alpha^{j}-1) \big\|_{\omega} \leq 3 \, \omega(j),
$$
ce qui entra\^\i ne (\ref{unite1}). Pour (\ref{unite2}), on pose $f_{m} = \sum \limits_{j=-m}^{+m} \hat{f}(j) (\alpha^{j} - 1)$, de sorte que $\dsp \lim_{m \rightarrow +\infty} \big\| f - f_{m} \big\|_{\omega} = 0$. On a alors en utilisant (\ref{unite1})
\begin{eqnarray} \label{unite11} 
\big\| (e_{n}-1)f \big\|_{\omega} & \leq & \big\| (e_{n}-1)f_{m} \big\|_{\omega} + \big\| (e_{n}-1)(f-f_{m}) \big\|_{\omega} \nonumber \\
& \leq & \big\| (e_{n}-1)f_{m} \big\|_{\omega} + 3 \, \big\| f-f_{m} \big\|_{\omega}. 
\end{eqnarray}
De (\ref{unite01}) et (\ref{unite02}), on d\'eduit que pour $j \in \bbz$, $\dsp \lim_{n \rightarrow +\infty} \big\| (e_{n}-1) (\alpha^{j}-1) \big\|_{\omega} = 0$, ce qui entra\^\i ne que $\dsp \lim_{n \rightarrow +\infty} \big\| (e_{n}-1)f_{m} \big\|_{\omega} = 0$. Il est maintenant facile de voir que $\dsp \lim_{n \rightarrow +\infty} \big\| (e_{n}-1)f \big\|_{\omega}  = 0$.
\end{proof}

\begin{prop} \label{ditkin}
Soit $\omega$ un poids satisfaisant $(W_{s})$, on pose
\begin{eqnarray} \label{unites}
u_{n} = e_{n}^{[s]+1}  \qquad (n \geq 1) ,
\end{eqnarray}
o\`u $e_{n}$ est la fonction d\'efinie en (\ref{unite0}). Alors pour toute fonction $f$ dans $A_{\omega}(\bbt)$ telle que $f^{(k)}(1)=0$ pour $0 \leq k \leq [s]$, on a
$$
\lim_{n \rightarrow +\infty} \big\| (u_{n}-1)f \big\|_{\omega} = 0.
$$ 
\end{prop}

\begin{proof}
Dans cette d\'emonstration on posera $p=[s]$ et on notera, pour $0 \leq k \leq p$, $\omega_{k} = \big( \omega_{k}(n) \big)_{n \in \bbz}$ la suite d\'efinie par
$$
\omega_{k}(n) = \frac{\omega(n)}{(1+|n|)^{k}}.
$$
En utilisant les relations $e_{n} = n (\alpha - 1)(e_{n} - 1)$ et $(\alpha-1) e_{n}' = - (e_{n} - 1)e_{n} \quad (n \geq 1)$, on montre par r\'ecurrence que pour tout $k$ dans $\{ 0, \ldots, p \}$, on a
\begin{eqnarray} \label{ditkin01}
(u_{n}-1)^{(k)} = \frac{(e_{n} - 1)}{(\alpha - 1)^{k}} P_{k} (e_{n}) \qquad (n \geq 1), 
\end{eqnarray}
o\`u $P_{k}$ est un polyn\^ome de degr\'e inf\'erieur \`a $k+p$ ne d\'ependant pas de $n$. \\
Soit maintenant $f$ dans $A_{\omega}(\bbt)$ telle que $f(1)=0$. D'apr\`es le lemme \ref{derivee}, il s'agit de montrer que 
\begin{eqnarray}
& & \lim_{n \rightarrow +\infty} \Big\| \big[ (u_{n}-1)f \big]^{(p)} \Big\|_{\omega_{p}} = 0 \label{ditkin02} \\
& & \lim_{n \rightarrow + \infty} \widehat{\big[ (u_{n}-1)f \big]}(k) = 0 \qquad \big( 0 \leq k \leq p-1 \big). \label{ditkin03}
\end{eqnarray}
Les conditions (\ref{ditkin03}) se montrent facilement \`a l'aide du th\'eor\`eme de convergence domin\'ee, il reste donc \`a montrer (\ref{ditkin02}). En utilisant la formule de Leibnitz et l'identit\'e (\ref{ditkin01}), on obtient
\begin{eqnarray}
\big[ (u_{n}-1)f \big]^{(p)} & = & \sum_{k=0}^{p} \binom{p}{k} (u_{n}-1)^{(k)} f^{(p-k)} \nonumber \\
& = & \sum_{k=0}^{p} \binom{p}{k} (e_{n} - 1) \frac{f^{(p-k)}}{(\alpha - 1)^{k}} P_{k} (e_{n}). \label{ditkin04} 
\end{eqnarray}
Or les hypoth\`eses sur le poids $\omega$ nous permettent d'utiliser le lemme \ref{derivee} qui nous assure que 
$$
\frac{f^{(p-k)}}{(\alpha - 1)^{k}} \in A_{\omega_{p}}(\bbt) \qquad (0 \leq k \leq p).
$$
De plus, $\omega_{p}$ v\'erifie la condition $(W_{\beta})$ avec $\beta=s-p$. On d\'eduit alors du lemme \ref{unite} que si $g \in A_{\omega_{p}}(\bbt)$ et $g(1)=0$, alors $\big\| e_{n} g \big\|_{\omega_{p}} \leq 4 \big\| g \big\|_{\omega_{p}}$ ($n \geq 1$). Et comme pour tout $k \in \{0, \ldots, p\}$, la fonction $\dsp \frac{f^{(p-k)}}{(\alpha - 1)^{k}}$ s'annule en $1$, on peut alors trouver une constante $C>0$ ind\'ependante de $n$ et de $k$ telle que
\begin{eqnarray} \label{ditkin05}
\Big\| (e_{n} - 1) \frac{f^{(p-k)}}{(\alpha - 1)^{k}} P_{k} (e_{n}) \Big\|_{\omega_{p}} \leq C \, \Big\| (e_{n} - 1) \frac{f^{(p-k)}}{(\alpha - 1)^{k}} \Big\|_{\omega_{p}}.
\end{eqnarray} 
De plus, toujours d'apr\`es le lemme \ref{unite},
\begin{eqnarray*}
\lim_{n \rightarrow +\infty} \Big\| (e_{n}-1) \frac{f^{(p-k)}}{(\alpha - 1)^{k}} \Big\|_{\omega_{p}} = 0 \qquad (0 \leq k \leq p). 
\end{eqnarray*}
Par cons\'equent, on d\'eduit alors de (\ref{ditkin04}) et (\ref{ditkin05}) que 
$$
\lim_{n \rightarrow +\infty} \Big\| \big[ (u_{n}-1)f \big]^{(p)} \Big\|_{\omega_{p}} = 0.
$$ 
\end{proof}

\section{Id\'eaux ferm\'es de $A_{\omega}(\bbt)$ et de $A_{s}^{+}(\bbt)$}

Il est montr\'e dans \cite{Atz1} que si $\omega$ est un poids v\'erifiant $(W_{s})$ et $(A)$, et $I$ un id\'eal ferm\'e de $A_{\omega}(\bbt)$ tel que $h^{0}(I) = \{ z_{0} \}$, alors il existe $j$ dans $\{ 1, \ldots, [s] \}$ tel que $I = \Big\{ f \in A_{\omega}(\bbt) : \, f(z_{0}) = \dots = f^{(j)}(z_{0}) = 0 \Big\}$. Nous allons montrer que ce r\'esultat s'\'etend aux ferm\'es d\'enombrables. Introduisons d'abord la notation suivante: si $f$ est dans $A_{\omega}(\bbt)$, et $I$ un id\'eal fem\'e de $A_{\omega}(\bbt)$, on pose $I(f) = \Big\{ g \in A_{\omega}(\bbt) : \, fg \in I \Big\}$. On a le r\'esultat suivant:

\begin{lemme}
Soit $\omega$ un poids v\'erifiant $(W_{s})$ et $(A)$, et $I$ un id\'eal ferm\'e de $A_{\omega}(\bbt)$ tel que $h^{0}(I)$ poss\`ede un point isol\'e $z_{0}$. Soit $k = \max \big\{ j \in \{0, \ldots, [s]\} : \, z_{0} \in h^{j}(I) \big\}$. Alors il existe $g$ dans $I$ de la forme $g = (\alpha - z_{0})^{k+1} \psi$, avec $\psi \in A_{\omega}(\bbt)$ et $\psi(z_{0}) \neq 0$.
\end{lemme}

\begin{proof}
$z_{0}$ \'etant un point isol\'e dans $h^{0}(I)$, il existe $\psi \in A_{\omega}(\bbt)$ telle que
$$
\psi = \left\{ \begin{array}{ll}
1 & \textrm{sur un voisinage de } z_{0} \\
0 & \textrm{sur un voisinage de } h^{0}(I) \backslash \{ z_{0} \}
\end{array} \right.. 
$$
L'alg\`ebre $A_{\omega}(\bbt)$ \'etant r\'eguli\`ere, $1-\psi \in I(\psi)$, ce qui prouve que $h^{0} \big( I(\psi) \big) \subset \{ z_{0} \}$. Puisque $\psi \notin I$, on a $I(\psi) \neq A_{\omega}(\bbt)$, et donc $h^{0} \big( I(\psi) \big) = \{ z_{0} \}$. Par cons\'equent, comme $\omega$ v\'erifie les conditions $(W_{s})$ et $(A)$, on d\'eduit de la proposition 6 de \cite{Atz1} qu'il existe un entier $\gamma$, $0 \leq \gamma \leq [s]$, tel que 
$$
I(\psi) = \Big\{ f \in A_{\omega}(\bbt) : \, f(z_{0}) = \dots = f^{(\gamma)}(z_{0}) = 0 \Big\}.
$$
Comme $I \subset I(\psi)$, on a $\gamma \leq k$. D'autre part $(\alpha - z_{0})^{\gamma+1}$ appartient \`a $I(\psi)$, et donc la fonction $g = (\alpha - z_{0})^{\gamma+1} \psi$ appartient \`a $I$. Comme $g^{(k)}(z_{0}) = 0$ et $\psi(z_{0}) \neq 0$, on a $\gamma \geq k$. Donc $\gamma = k$, et la fonction $g$ ainsi d\'efinie convient.
\end{proof}

\begin{theo} \label{ideauxdeaw}
Soit $\omega$ un poids v\'erifiant $(W_{s})$ et $(A)$, et $I$ un id\'eal ferm\'e de $A_{\omega}(\bbt)$ tel que $h^{0}(I)$ est d\'enombrable. Alors  
\begin{eqnarray*} 
I = \Big\{ f \in A_{\omega}(\bbt) : \, f^{(j)}=0 \textrm{ sur } h^{j}(I) \quad (0 \leq j \leq [s]) \Big\}. 
\end{eqnarray*}
\end{theo}

\begin{proof}
Soit $I$ un id\'eal ferm\'e de $A_{\omega}(\bbt)$ tel que $h^{0}(I)$ est d\'enombrable, notons $J = \Big\{ f \in A_{\omega}(\bbt) : \, f^{(j)}=0 \textrm{ sur } h^{j}(I) \quad (0 \leq j \leq [s]) \Big\}$. L'inclusion $I \subset J$ \'etant \'evidente, il reste \`a montrer l'autre. Soit $f \in J$, nous allons montrer que $I(f) = A_{\omega}(\bbt)$. Soit $z_{0} \in h^{0}(I) \backslash h^{[s]}(I)$ et $k \in \{0, \ldots, [s]-1\}$ tel que $z_{0} \in h^{k}(I) \backslash h^{k+1}(I)$. Il est facile de voir que $J \subset \overline{(\alpha - z_{0})^{k+1} A_{\omega}(\bbt)}^{\, \| \, \|_{\omega}}$, et donc il existe une suite $(f_{m})_{m \geq 0}$ de fonctions de la forme $f_{m} = (\alpha - z_{0})^{k+1} \phi_{m}$, avec $\phi_{m} \in A_{\omega}(\bbt)$, telle que $\dsp \lim_{m \rightarrow +\infty} \big\| f - f_{m} \big\|_{\omega} = 0$. De plus, comme $z_{0} \notin h^{[s]}(I)$, le point $z_{0}$ est n\'ecessairement isol\'e dans $h^{0}(I)$. Donc d'apr\`es le lemme pr\'ec\'edent, il existe $g$ dans $I$ qui s'\'ecrit $g = (\alpha - z_{0})^{k+1} \psi$ avec $\psi \in A_{\omega}(\bbt)$ et $\psi(z_{0}) \neq 0$. Posons alors 
$$
\Psi_{m} = \phi_{m} g = f_{m} \psi \quad (m \geq 0).
$$
On a pour tout entier $m \geq 0$, $\Psi_{m} \in I$ et $\dsp \lim_{m \rightarrow  +\infty} \big\| \Psi_{m} - f \psi \big\|_{\omega} = 0$. Comme $I$ est ferm\'e, on a donc $f \psi \in I$, c'est-\`a-dire $\psi \in I(f)$. Par cons\'equent, $z_{0} \notin h^{0} \big( I(f) \big)$ (car $\psi(z_{0}) \neq 0$). Finalement, on en d\'eduit dans un premier temps que 
\begin{eqnarray} \label{ideauxdeaw00}
h^{0} \big( I(f) \big) \subset h^{[s]}(I).
\end{eqnarray}      
Supposons que $h^{0} \big( I(f) \big) \neq \emptyset$, alors $h^{0} \big( I(f) \big)$ admet un point isol\'e $\xi_{0}$. Sans perte de g\'en\'eralit\'e, on va supposer que $\xi_{0} = 1$. L'alg\`ebre $A_{\omega}(\bbt)$ \'etant r\'eguli\`ere, il existe une fonction $\Phi$ telle que 
$$
\Phi = \left\{ \begin{array}{ll}
1 & \textrm{sur un voisinage de } 1 \\
0 & \textrm{sur un voisinage de } h^{0} \big( I(f) \big) \backslash \{ 1 \}
\end{array} \right.. 
$$
On pose $L_{\Phi} = \big\{ h \in A_{\omega}(\bbt) : \, h \Phi \in I(f) \big\}$, l'id\'eal de division de $I(f)$ par $\Phi$. On a $1-\Phi \in L_{\Phi}$, et donc $h^{0} \big( L_{\Phi} \big) \subset \{ 1 \}$. Comme $\omega$ v\'erifie les conditions $(A)$ et $(W_{s})$, on d\'eduit de la proposition 6 de \cite{Atz1} que $\Big\{ f \in A_{\omega}(\bbt) : \, f(1) = \dots = f^{(\gamma)}(1) = 0 \Big\} \subset L_{\Phi}$. Par cons\'equent la suite $\big( u_{n} \big)_{n \geq 1}$ d\'efinie en (\ref{unites}) appartient \`a $L_{\Phi}$. D'apr\`es l'inclusion $h^{0} \big( I(f) \big) \subset h^{[s]}(I)$ \'etablie en (\ref{ideauxdeaw00}), on a $f^{(k)}(1)=0$ pour $0 \leq k \leq [s]$. Et comme le poids $\omega$ v\'erifie la condition $(W_{s})$, on d\'eduit du lemme \ref{ditkin} que $\dsp \lim_{n \rightarrow +\infty} \big\| (u_{n}-1)f \big\|_{\omega} = 0$. Puisque $u_{n} \Phi f \in I \,\, (n \geq 1)$ et $\dsp \lim_{n \rightarrow +\infty} \big\| u_{n} \Phi f - \Phi f \big\|_{\omega} = 0$, on a $\Phi f \in I$, ce qui contredit le fait que $1 \in h^{0} \big( I(f) \big)$. Finalement on a montr\'e que $h^{0} \big( I(f) \big) = \emptyset$, et donc $I(f) = A_{\omega}(\bbt)$, ce qui signifie exactement que $f \in I$.  
\end{proof}

On va maintenant s'int\'eresser aux id\'eaux ferm\'es de $A_{s}^{+}(\bbt)$. Puisque l'alg\`ebre $A_{s}(\bbt)$ v\'erifie la condition de Ditkin analytique forte, le th\'eor\`eme B de \cite{BeGi} nous assure que tout id\'eal ferm\'e $I$ de $A_{s}^{+}(\bbt)$ tel que $h_{+}^{0}(I)$ est fini, est de la forme 
$$
I = I^{A_{s}} \cap S_{I} \calh^{\infty}(\bbd),
$$
o\`u $I^{A_{s}}$ est l'id\'eal ferm\'e de $A_{s}(\bbt)$ engendr\'e par $I$. Par cons\'equent, compte-tenu du th\'eor\`eme \ref{ideauxdeaw}, $I$ est de la forme 
\begin{eqnarray*} 
I = \Big\{ f \in A_{s}^{+}(\bbt) : \, S_{I} \big| S(f) \textrm{ et } f^{(j)}=0 \textrm{ sur } h^{j}(I) \cap \bbt \quad (0 \leq j \leq [s]) \Big\}. 
\end{eqnarray*}

Soit $I$ un id\'eal ferm\'e non nul de $A_{s}^{+}(\mathbb{T})$, on note
$$
\pi_{s}^{+} : A_{s}^{+}(\mathbb{T}) \longrightarrow A_{s}^{+}(\mathbb{T})/I
$$
la surjection canonique. On a alors le r\'esultat suivant:

\begin{lemme} \label{evaluer}
Soient $s$ un r\'eel positif et $I$ un id\'eal ferm\'e non r\'eduit \`a $\{ 0 \}$ de $A_{s}^{+}(\mathbb{T})$ tel que $S_{I}=1$. Alors
$$ 
\big\| \pi_{s}^{+}(\alpha)^{-n} \big\| = O \big( e^{ \varepsilon \sqrt{n}} \big) \,\, (n \rightarrow +\infty), \, \textrm{ pour tout } \varepsilon > 0. \\    
$$
\end{lemme}

\begin{proof}
C'est un r\'esultat \'etabli par A. Atzmon dans la preuve de la proposition 8 de \cite{Atz1} dans le cas o\`u $I = \Big\{ f \in A_{s}^{+}(\mathbb{T}) : \, f_{|_{E}} = 0 \Big\}$, et qui est une cons\'equence du lemme 5.c de \cite{Atz1} (voir aussi \cite{Zar1}, proposition 2.1). Le r\'esultat ci-dessus se d\'emontre de fa\c con analogue. 
\end{proof}

Nous avons alors le th\'eor\`eme suivant:

\begin{theo} \label{caracden}
Soient $s$ un r\'eel positif et $I$ un id\'eal ferm\'e de $A_{s}^{+}(\bbt)$ tel que $S_{I}=1$ et $h_{+}^{0}(I)$ est au plus d\'enombrable. Alors  
\begin{eqnarray} \label{lacarac}
I = I \big( 1 ; h_{+}^{0}(I), \ldots, h_{+}^{[s]}(I) \big). 
\end{eqnarray} 
\end{theo}

\begin{proof}
Soit $I$ un id\'eal ferm\'e sans facteur int\'erieur de $A_{s}^{+}(\mathbb{T})$ tel que $h_{+}^{0}(I)$ est au plus d\'enombrable. On d\'eduit du lemme \ref{evaluer} que
$$
\big\| \pi_{s}^{+}(\alpha)^{-n} \big\| = O \big( e^{\varepsilon \sqrt{n}} \big) \, (n \rightarrow +\infty), \, \textrm{ pour tout } \varepsilon >  0.
$$
On consid\`ere alors le poids $\omega$ d\'efini par 
$$
\left\{ \begin{array}{ll}
\omega(n) = (1 + n)^{s} & (n \geq 0) \\
\omega(-n) = (1 + n)^{s} \displaystyle \sup_{0 < k \leq n} \big\| \pi_{s}^{+}(\alpha)^{-k} \big\| & (n > 0),
\end{array} \right.
$$
et on d\'efinit l'application continue $\theta : A_{\omega}(\mathbb{T}) \longrightarrow A_{s}^{+}(\mathbb{T})/I$ par
$$
\theta (f) = \sum_{n = -\infty}^{+\infty} \hat{f}(n) \pi_{s}^{+}(\alpha)^{n}.
$$
On a $\theta_{|_{A_{s}^{+}(\mathbb{T})}} = \pi_{s}^{+}$, et donc 
\begin{eqnarray} \label{sansfac1}
\textrm{Ker}\, \theta \cap A_{s}^{+}(\mathbb{T}) = I.
\end{eqnarray}
Si $I^{A_{\omega}}$ d\'esigne l'id\'eal ferm\'e de $A_{\omega}(\bbt)$ engendr\'e par $I$, on a $\dsp I^{A_{\omega}} = \overline{\bigcup_{n \geq 0} \alpha^{-n} I}^{A_{\omega}(\bbt)}$. Or il est facile de voir que pour tout $n \geq 0$, $\alpha^{-n} I \subset \textrm{Ker}\,\theta$, et donc 
$$
I^{A_{\omega}} \cap A_{s}^{+}(\bbt) \subset \textrm{Ker}\,\theta \cap A_{s}^{+}(\bbt) = I.
$$
L'autre inclusion \'etant \'evidente, on a donc
$$
I = I^{A_{\omega}} \cap A_{s}^{+}(\bbt).
$$
Puisque $S_{I}=1$, on a pour tout $0 \leq k \leq [s]$, $h_{+}^{k}(I) \subset \bbt$. Il est alors facile de voir que pour tout $0 \leq k \leq [s]$, $h^{k} \big( I^{A_{\omega}} \big) = h_{+}^{k}(I)$, on d\'eduit finalement du th\'eor\`eme \ref{ideauxdeaw} que
$$
I = I \big( 1 ; h_{+}^{0}(I), \ldots, h_{+}^{[s]}(I) \big). 
$$ 
\end{proof}

\section{Interpolation}

On note $\calc^{\infty}(\bbt)$ l'ensemble des fonctions de classe $\calc^{\infty}$ sur le cercle, que l'on \'equipe de sa topologie d'espace de Fr\'echet usuelle d\'efinie par les normes $\big( \rho_{\nu} \big)_{\nu \in \bbr^{+}}$ d\'efinies par
$$
\rho_{\nu}(f) = \sum_{n=-\infty}^{+\infty} \big| \hat{f}(n) \big| (1+|n|)^{\nu}.
$$
On note \'egalement
$$
\textrm{{\LARGE $a$}}^{\infty} = \textrm{{\LARGE $a$}}^{\infty}(\bbd) = \Big\{ f \textrm{ holomorphe dans } \bbd : \, f \in \calc(\overline{\bbd}) \textrm{ et } f_{|_{\bbt}} \in \calc^{\infty}(\bbt) \Big\},
$$
que l'on regardera comme l'ensemble des fonctions de $\calc^{\infty}(\bbt)$ \`a coefficients de Fourier strictement n\'egatifs nuls. Soit $E$ un ferm\'e du cercle unit\'e, on d\'efinit 
$$
I_{\infty}(E) = \Big\{ f \in \calc^{\infty}(\bbt) : \, f_{|_{E}}^{(i)} = 0 \quad (i \geq 0) \Big\}, 
$$
et on pose $I_{\infty}^{+}(E) = I_{\infty}(E) \cap \textrm{{\LARGE $a$}}^{\infty}$. \\
Le dual de $\calc^{\infty}(\bbt)$ est $\mathcal{D}'(\bbt)$, l'ensemble des distributions sur $\bbt$. On associe \`a chaque distribution $T \in \mathcal{D}'(\bbt)$ une suite de coefficients de Fourier $\big( \widehat{T}(n) \big)_{n \in \bbz}$, o\`u pour tout entier $n$, $\hat{T}(n) = \big< \alpha^{-n},T \big>$ (avec $\alpha: z \mapsto z$), qui v\'erifient $\big| \hat{T}(n) \big| = O(|n|^{m})$ pour un certain entier $m \geq 0$. La dualit\'e entre $\calc^{\infty}(\bbt)$ et $\mathcal{D}'(\bbt)$ est donn\'ee par la formule
$$
\big< f , T\big> = \sum_{n = -\infty}^{+\infty} \hat{f}(n) \widehat{T}(-n) \qquad \big( f \in \calc^{\infty}(\bbt) , T \in \mathcal{D}'(\bbt) \big). 
$$
$I_{\infty}(E)^{\perp}$ \big(resp. $I_{\infty}^{+}(E)^{\perp}$\big) d\'esigne l'ensemble des distributions s'annulant sur $I_{\infty}(E)$ \big(resp. sur $I_{\infty}^{+}(E)$\big). De m\^eme $\big( \textrm{{\LARGE $a$}}^{\infty} \big)^{\perp}$ est l'ensemble des distributions s'annulant sur $\textrm{{\LARGE $a$}}^{\infty}$, c'est-\`a-dire les distributions \`a coefficients n\'egatifs nuls. \\

Dans le cas particulier o\`u le poids $\omega$ est d\'efini par $\omega(n) = (1+n)^{s}$ et $\omega(-n) = (1+n)^{t}$ pour tout $n \geq 0$, avec $t$ et $s$ deux r\'eels, on notera $\big( A_{s,t}(\bbt), \| \, \|_{s,t} \big)$ l'alg\`ebre $\big( A_{\omega}(\bbt), \| \, \|_{\omega} \big)$. Nous supposerons dor\'enavant que $t \geq s$, de sorte que $A_{s,t}(\bbt)$ satisfait les conditions $(W_{s})$ et $(A)$. Si $E$ est un ferm\'e de $\bbt$, on notera
$$ 
I_{s,t}(E) = \Big\{ f \in A_{s,t}(\bbt) : \, f_{|_{E}} = \ldots = f_{|_{E}}^{([s])} = 0 \Big\},
$$
et $I_{s}^{+}(E) = I_{s,t}(E) \cap A_{s}^{+}(\bbt)$. On identifiera le dual de $A_{s,t}(\bbt)$ (que l'on note $\big( A_{s,t}(\bbt) \big)'$) au sous-espace de $\mathcal{D}'(\bbt)$ form\'e des distributions $T$ telles que $\dsp \sup_{n \leq 0} \frac{|\hat{T}(n)|}{(1+|n|)^{s}} + \sup_{n > 0} \frac{|\hat{T}(n)|}{(1+n)^{t}} < +\infty$. $I_{s,t}(E)^{\perp}$ \big(resp. $I_{s}^{+}(E)^{\perp}$\big) d\'esigne l'ensemble des \'el\'ements de $\big( A_{s,t}(\bbt) \big)'$ s'annulant sur $I_{s,t}(E)$ \big(resp. sur $I_{s}^{+}(E)$\big). De m\^eme $A_{s}^{+}(\bbt)^{\perp}$ est l'ensemble des \'el\'ements de $\big( A_{s,t}(\bbt) \big)'$ s'annulant sur $A_{s}^{+}(\bbt)$, c'est-\`a-dire les \'el\'ements de $\big( A_{s,t}(\bbt) \big)'$ \`a coefficients n\'egatifs nuls. \\

\begin{defi} \label{defiinter}
Soient $s$ et $t$ deux r\'eels positifs tels que $t \geq s$. On dira qu'un ferm\'e $E$ du cercle unit\'e est d'interpolation pour $A_{s,t}(\bbt)$ si 
$$
\forall f \in A_{s,t}(\bbt) , \quad \exists \, g \in A_{s}^{+} (\bbt) \,: \quad f_{|_{E}}^{(i)} = g_{|_{E}}^{(i)} \qquad (0 \leq i \leq [s]),
$$
c'est-\`a-dire si $A_{s}^{+} (\bbt) + I_{s,t}(E) = A_{s,t}(\bbt)$. \\
On dira qu'un ferm\'e $E$ du cercle unit\'e est d'interpolation pour $\textrm{{\LARGE $a$}}^{\infty}$ si 
$$
\forall f \in \calc^{\infty}(\bbt) , \quad \exists \, g \in \textrm{{\LARGE $a$}}^{\infty} \, : \quad f_{|_{E}}^{(i)} = g_{|_{E}}^{(i)} \qquad (i \geq 0),
$$
c'est-\`a-dire si $\textrm{{\LARGE $a$}}^{\infty} + I_{\infty}(E) = \calc^{\infty}(\bbt)$.
\end{defi}

H. Alexander, B. A. Taylor et D. L. Williams ont donn\'e une caract\'erisation g\'eom\'etrique des ensembles d'interpolation pour $\textrm{{\LARGE $a$}}^{\infty}$ (\cite{ATaWi2}). Ils ont montr\'e qu'un ferm\'e du cercle unit\'e est d'interpolation pour $\textrm{{\LARGE $a$}}^{\infty}$ si et seulement si $E$ v\'erifie la condition (ATW).
On a \'egalement la caract\'erisation suivante de ces ensembles:

\begin{prop} \label{interc}
Soit $E$ un ferm\'e du cercle unit\'e, alors les propri\'et\'es suivantes sont \'equivalentes: \\
\hspace*{0.3cm} (i) $E$ est un ensemble d'interpolation pour $\textrm{{\LARGE $a$}}^{\infty}$. \\
\hspace*{0.2cm} (ii) $I_{\infty}^{+}(E)^{\perp} = I_{\infty}(E)^{\perp} + \big( \textrm{{\LARGE $a$}}^{\infty} \big)^{\perp}$. \\
\hspace*{0.1cm} (iii) Pour tout $s \geq 0$, il existe une constante $C > 0$ et $t \geq 0$ tels que 
$$
\sup_{n \in \bbz} \frac{|\hat{T}(n)|}{(1+|n|)^{t}} \leq C \sup_{n \leq 0} \frac{|\hat{T}(n)|}{(1+|n|)^{s}} \qquad \Big( T \in I_{\infty}(E)^{\perp} \Big).
$$
\end{prop}

\begin{proof}
L'\'equivalence $(i) \! \iff \! (ii)$ a \'et\'e \'etablie dans \cite{ATaWi2} (proposition 2.1). On va achever la preuve en prouvant que $(i) \! \iff \! (iii)$.
On sait que $E$ est d'interpolation pour $\textrm{{\LARGE $a$}}^{\infty}$ si et seulement si l'injection canonique 
$$
i \, : \, \textrm{{\LARGE $a$}}^{\infty} / I_{\infty}^{+}(E) \longrightarrow \calc^{\infty}(\bbt) / I_{\infty}(E)
$$
est surjective. Puis en utilisant la proposition 4, IV.30 de \cite{BourEVT4} qui caract\'erise les surjections entre espaces de Fr\'echet, on d\'eduit que E est d'interpolation pour $\textrm{{\LARGE $a$}}^{\infty}$ si et seulement si pour tout $s \geq 0$, il existe une constante $C >0$ et $t \geq 0$ tels que pour tout $T$ dans $I_{\infty}(E)^{\perp}$
\begin{eqnarray*}
\Big( \big| \big< f,T \big> \big| \leq \rho_{s}(f) \, , \quad f \in \textrm{{\LARGE $a$}}^{\infty} / I_{\infty}^{+}(E) \Big) & \Longrightarrow & \Big( \big| \big< f,T \big> \big| \leq C \rho_{t}(f) \, , \quad f \in \calc^{\infty}(\bbt) / I_{\infty}(E) \Big),
\end{eqnarray*}
c'est-\`a-dire si et seulement si pour tout $s \geq 0$, il existe une constante $C >0$ et $t \geq 0$ tels que 
\begin{eqnarray*}
\Big( \sup_{n \leq 0} \frac{|\hat{T}(n)|}{(1+|n|)^{s}} \leq 1 \, , \quad T \in I_{\infty}(E)^{\perp} \Big) & \Longrightarrow & \sup_{n \in \bbz} \frac{|\hat{T}(n)|}{(1+|n|)^{t}} \leq C,
\end{eqnarray*}
ce qui est clairement \'equivalent \`a $(iii)$.
\end{proof}

On a \'egalement un r\'esultat analogue concernant les ensembles d'interpolation pour $A_{s,t}(\bbt)$ ($t \geq s$):

\begin{prop} \label{intera}
Soit $E$ un ferm\'e du cercle unit\'e et $s,t$ deux r\'eels positifs tels que $t \geq s$. Alors les propri\'et\'es suivantes sont \'equivalentes: \\
\hspace*{0.3cm} (i) $E$ est un ensemble d'interpolation pour $A_{s,t}(\bbt)$. \\
\hspace*{0.2cm} (ii) $I_{s}^{+}(E)^{\perp} = I_{s,t}(E)^{\perp} + A_{s}^{+}(\bbt)^{\perp}$. \\
\hspace*{0.1cm} (iii) Il existe une constante $C > 0$ telle que 
$$
\sup_{n > 0} \frac{|\hat{T}(n)|}{(1+|n|)^{t}} + \sup_{n \leq 0} \frac{|\hat{T}(n)|}{(1+|n|)^{s}} \leq C \sup_{n \leq 0} \frac{|\hat{T}(n)|}{(1+|n|)^{s}} \qquad \Big( T \in I_{s,t}(E)^{\perp} \Big).
$$
\end{prop}

On est maintenant en mesure d'\'enoncer un r\'esultat qui fait le lien entre les ensembles d'interpolation pour $\textrm{{\LARGE $a$}}^{\infty}$ et les ensembles d'interpolation pour $A_{s,t}(\bbt)$ ($t \geq s$).
 
\begin{theo} \label{interdeux}
Soit $E$ un ferm\'e du cercle unit\'e $\bbt$. \\
\hspace*{0.2cm} 1) On suppose que $E$ est d'interpolation pour $\textrm{{\LARGE $a$}}^{\infty}$. Alors pour tout $s \geq 0$, il existe $t \geq s$ tel que $E$ soit d'interpolation pour $A_{s,t}(\bbt)$. \\
\hspace*{0.2cm} 2) R\'eciproquement, on suppose que $E$ est d\'enombrable et que pour tout $s \geq 0$, il existe $t \geq s$ tel que $E$ soit d'interpolation pour $A_{s,t}(\bbt)$. Alors $E$ est d'interpolation pour $\textrm{{\LARGE $a$}}^{\infty}$.
\end{theo}

\begin{proof}
1) Cette assertion d\'ecoule directement des caract\'erisations $(iii)$ des propositions ($\ref{interc}$) et ($\ref{intera}$) puisque $I_{s,t}(E)^{\perp}$ s'injecte dans $I_{\infty}(E)^{\perp}$. \\
2) Supposons que $E$ v\'erifie les hypoth\`eses de l'assertion 2), on va montrer que $E$ satisfait la propri\'et\'e $(iii)$ de la proposition \ref{interc}. Soit $s$ un r\'eel positif, il existe $t \geq s$ tel que $E$ soit d'interpolation pour $A_{s,t}(\bbt)$. Donc d'apr\`es la proposition \ref{intera}, il existe une constante $C > 0$ telle que pour tout \'el\'ement de $I_{s,t}(E)^{\perp}$,
\begin{eqnarray} \label{interdeux1}
\sup_{n > 0} \frac{|\hat{T}(n)|}{(1+|n|)^{t}} \leq C \sup_{n \leq 0} \frac{|\hat{T}(n)|}{(1+|n|)^{s}}.
\end{eqnarray}
Soit alors $T$ un \'el\'ement de $I_{\infty}(E)^{\perp}$. Si $\dsp \sup_{n \leq 0} \frac{|\hat{T}(n)|}{(1+|n|)^{s}} = +\infty$, alors $T$ v\'erifie trivialement (\ref{interdeux1}). Supposons donc que $\dsp \sup_{n \leq 0} \frac{|\hat{T}(n)|}{(1+|n|)^{s}} < +\infty$. Il existe alors $t' \geq t$ tel que $T$ soit dans $\big( A_{s,t'}(\bbt) \big)'$. Notons $K$ la fermeture de ${I_{\infty}(E)}$ dans $A_{s,t'}(\bbt)$. On v\'erifie facilement que $K$ est un id\'eal ferm\'e de $A_{s,t'}(\bbt)$ tel que $h(K) \subset E$, et ainsi d'apr\`es le th\'eor\`eme (\ref{ideauxdeaw}), $I_{s,t'}(E) \subset K$. L'inclusion inverse \'etant acquise, on a finalement $K = I_{s,t'}(E)$. Et comme $T$ est dans $I_{\infty}(E)^{\perp}$, $T$ appartient \`a $I_{s,t'}(E)^{\perp}$ par continuit\'e de $T$ sur $A_{s,t'}(\bbt)$. Mais puisque $E$ est d'interpolation pour $A_{s,t}(\bbt)$, les inclusions naturelles
$$
A_{s}^{+}(\bbt) / I_{s}^{+}(E) \subset A_{s,t'}(\bbt) / I_{s,t'}(E) \subset A_{s,t}(\bbt) / I_{s,t}(E),
$$
sont en r\'ealit\'e des \'egalit\'es. Par cons\'equent leurs duaux sont \'egaux et donc $T \in I_{s,t}(E)^{\perp}$. Ainsi $T$ v\'erifie la propri\'et\'e (\ref{interdeux1}), ce qui ach\`eve la d\'emonstration.
\end{proof}

\section{Interpolation dans $A_{s,t}(\bbt \, )$ et op\'erateurs}

\begin{prop} \label{interop}
Soit $E$ un ferm\'e du cercle unit\'e, $s$ et $t$ deux r\'eels positifs tels que $t \geq s$. \\
\hspace*{0.2cm} 1) On suppose que $E$ est d\'enombrable et d'interpolation pour $A_{s,t}(\bbt)$, alors $E$ v\'erifie la propri\'et\'e $P(s,t,E)$. \\
\hspace*{0.2cm} 2) R\'eciproquement, si $E$ v\'erifie la condition de Carleson (C) et la propri\'et\'e $P(s,t,E)$, alors $E$ est d'interpolation pour $A_{s,t}(\bbt)$.
\end{prop}

\begin{proof}
1) Soit $T$ un op\'erateur inversible sur un espace de Banach $X$ v\'erifiant les condition (\ref{introCp}) et (\ref{introA}). On d\'efinit alors le poids $\omega$ par 
$$
\left\{ \begin{array}{ll}
\omega(n) = (1 + n)^{s} & (n \geq 0) \\
\omega(-n) = (1 + n)^{s} \dsp \sup_{0 < k \leq n} \big\| T^{-k} \big\| & (n > 0).
\end{array} \right.
$$
Puisque $T$ v\'erifie la condition (\ref{introCp}), on peut d\'efinir un op\'erateur born\'e $\Phi$ de $A_{\omega}(\bbt)$ dans $\mathcal{L}(X)$ par 
$$
\Phi(f) = f(T) = \sum_{n = - \infty}^{+\infty} \hat{f}(n) T^{n} \qquad \big( f \in A_{\omega}(\bbt) \big).
$$
Puisque $A_{\omega}(\bbt)$ est r\'eguli\`ere, on a $h(\textrm{Ker} \, \Phi) = \textrm{Sp} \, T \subset E$ (voir \cite{EStZo1}, th\'eor\`eme 2.5), et donc $ \Big\{ f \in A_{\omega}(\bbt) : \, f_{|_{E}} = \dots = f^{([s])}_{|_{E}} = 0 \Big\} \subset \textrm{Ker} \, \Phi$ d'apr\`es le th\'eor\`eme \ref{ideauxdeaw}. En utilisant le fait que $E$ est d'interpolation pour $A_{s,t}(\bbt)$ et des m\'ethodes similaires \`a celles utilis\'ees dans \cite{Zar2} pour la preuve du th\'eor\`eme 2.6, on montre que, pour $n \geq 0$,
\begin{eqnarray*}
\big\| T^{-n} \big\| \leq C \big\| \alpha^{-n} \big\|_{s,t} = C (1+n)^{t},
\end{eqnarray*}
ce qui prouve que $T$ v\'erifie la propri\'et\'e (\ref{introDq}). \\
2) Pour la r\'eciproque, on consid\`ere l'op\'erateur $T$ d\'efini sur $A_{s}^{+}(\bbt)/I_{s}^{+}(E)$ par
$$ 
T \, : \, \pi_{s}^{+}(f) \longmapsto \pi_{s}^{+}(\alpha f) \qquad \big( f \in A_{s}^{+}(\bbt) \big), 
$$
o\`u $\pi_{s}^{+}$ est la surjection canonique de $A_{s}^{+}(\bbt)$ sur $A_{s}^{+}(\bbt)/I_{s}^{+}(E)$. On a $\big\| T^{n} \big\| = O \big( n^{s} \big)$ ($n \rightarrow +\infty$). Or, puisque $E$ est de Carleson, $I_{s}^{+}(E)$ est non r\'eduit \`a $\{ 0 \}$ (voir \cite{Carl}). De plus, $I_{s}^{+}(E)$ est sans facteur int\'erieur, donc en utilisant le lemme \ref{evaluer}, on obtient l'\'evaluation  
\begin{eqnarray} \label{evaluation}
\big\| T^{-n} \big\| = O \big( e^{\e \sqrt{n}} \big) \,\, (n \rightarrow +\infty), \, \textrm{ pour tout } \e > 0.    
\end{eqnarray}
Et puisque $\textrm{Sp} \, T = \textrm{Sp} \, \pi_{s}^{+}(\alpha) = E$, on a
$$ 
\big\| T^{-n} \big\| = \big\| \pi_{s}^{+}(\alpha)^{-n} \big\| = O \big( n^{t} \big) \, (n \rightarrow +\infty), 
$$
ce qui montre que $E$ est un ensemble d'interpolation pour $A_{s,t}(\bbt)$.
\end{proof}

Le r\'esultat suivant, annonc\'e dans l'introduction, est alors une cons\'equence directe du th\'eor\`eme \ref{interdeux} et de la proposition \ref{interop}. 

\begin{theo} \label{thetheo}
Soit $E$ un ferm\'e d\'enombrable du cercle unit\'e $\bbt$, alors les conditions suivantes sont \'equivalentes: \\
\hspace*{3mm} (i) $E$ v\'erifie la condition (ATW). \\
\hspace*{2mm} (ii) $E$ v\'erifie la condition (C) et pour tout $s \geq 0$, il existe $t \geq s$ tel que la propri\'et\'e $P(s,t,E)$ soit v\'erifi\'ee.
\end{theo}

Nous allons conclure en montrant que les hypoth\`eses du th\'eor\`eme sont optimales. Soit $E$ un ensemble ferm\'e du cercle unit\'e et $\mu$ une mesure \`a support dans $E$. Soit $J_{\mu}$ la fonction singuli\`ere associ\'ee \`a $\mu$, \`a savoir 
$$
J_{\mu}(z) = \exp \Big\{ -\frac{1}{2 \pi} \int_{0}^{2 \pi} \frac{e^{it} + z}{e^{it} - z} \mathrm{d} \mu(t) \Big\} \qquad (|z| < 1).
$$
On pose $\calh_{0} = \calh^{2}(\bbd) \ominus J_{\mu} \calh^{2}(\bbd)$, o\`u $\calh^{2}(\bbd)$ d\'esigne l'espace de Hardy usuel. On note $P_{\calh_{0}}$ la projection orthogonale sur $\calh_{0}$ et $\alpha \, : z \mapsto z$. On d\'efinit alors l'op\'erateur $T_{\mu}$ sur $\calh_{0}$ par
\begin{eqnarray} \label{optimal1}
T_{\mu}(f) = P_{\calh_{0}}(\alpha f) \qquad (f \in \calh_{0}).
\end{eqnarray}
D'apr\`es la proposition 5.1 p. 117 de \cite{NaFo}, on a $\textrm{Sp} \, T_{\mu} = \textrm{Supp} \, \mu \subset E$. \\
Si on remplace la condition ''pour tout $\e > 0$'' dans la propri\'et\'e $P(s,t,E)$ par une condition \`a $\e$ fix\'e, alors la propri\'et\'e $P(0,t,\{ 1 \})$ n'est v\'erifi\'ee pour aucun r\'eel $t \geq 0$. En effet, soient $\e_{0} > 0$, $\mu = 2 \pi \e_{0}^{2} \, \delta_{1}$ (o\`u $\delta_{1}$ est la mesure de Dirac en $1$) et $T_{1} = T_{\mu}$ l'op\'erateur d\'efini en (\ref{optimal1}). $T_{1}$ est un op\'erateur non unitaire tel que $\textrm{Sp} \, T_{1} = \{1\}$ et 
$$
\big\| T_{1}^{-n} \big\| = O(e^{4 \e_{0} \sqrt{n}}) \,\, (n \rightarrow +\infty) 
$$
(voir \cite{Zar2} pour plus de d\'etails). Le th\'eor\`eme 6.4 de \cite{Est1} nous montre alors que $T_{1}^{-n} = O(n^{t}) \, (n \rightarrow +\infty)$ n'est satisfait pour aucun r\'eel $t \geq 0$. \\
Nous allons maintenant montrer que l'hypoth\`ese de d\'enombrabilit\'e de $E$ dans le th\'eor\`eme \ref{thetheo} est essentielle. Pour cela, on a besoin du lemme suivant:

\begin{lemme} \label{contcarl}
Tout ferm\'e non d\'enombrable et de mesure nulle du cercle unit\'e contient un ensemble parfait qui v\'erifie la condition de Carleson (C).
\end{lemme}

\begin{proof}
Soit $S$ un ferm\'e non d\'enombrable et de mesure nulle du cercle unit\'e. Sans perte de g\'eneralit\'e, on peut supposer que $1 \notin S$. On \'ecrit alors $S = \big\{ e^{it} : \, t \in E \big\}$, o\`u $E$ est un ferm\'e de $\bbr$ inclu dans $]0, 2 \pi[$, et on se ram\`ene ainsi \`a la droite r\'eelle. Soit $P$ la partie parfaite de $E$. On pose
$$
P_{0} = \Big\{ x \in P : \, \exists \, \e > 0 \, \textrm{ tel que } \, ] x - \e, x [ \, \cap P = \emptyset \textrm{ ou } ] x, x + \e [ \, \cap P = \emptyset \Big\},
$$
c'est-\`a-dire l'ensemble des points de $P$ qui sont limites d'un seul c\^ot\'e d'une suite de points de $P$. Montrons dans un premier temps que $P_{0}$ est au plus d\'enombrable. On pose pour cela
$$
Q_{n}  = \Big\{ x \in P_{0} : \, \big] x - \frac{1}{n}, x \big[ \, \cap P = \emptyset \textrm{ ou } \big] x, x + \frac{1}{n} \big[ \, \cap P = \emptyset \Big\} \quad (n \geq 1).
$$
Il est clair que $P_{0} = \bigcup \limits_{n=1}^{+\infty} Q_{n}$ et que chaque $Q_{n}$ est fini. Par cons\'equent $P_{0}$ est au plus d\'enombrable. \\
Soient $a_{0}$ et $b_{0}$ deux points de $P \backslash P_{0}$ tels que $a_{0} < b_{0}$ et $\big| a_{0} - b_{0} \big| \leq 1$, on pose $I_{0} = [a_{0},b_{0}]$. A la premi\`ere \'etape, on retire de $I_{0}$ un intervalle ouvert $J_{1}^{(1)}$ dont les extr\'emit\'es sont dans $P \backslash P_{0}$, de sorte qu'il reste deux intervalles ferm\'es $I_{1}^{(1)}$ et $I_{2}^{(1)}$ qui soient non vides, non r\'eduits \`a un singleton et de longueur inf\'erieure \`a $\dsp \frac{1}{3} (b_{0} - a_{0})$ (un tel choix est possible car $P_{0}$ est au plus d\'enombrable). On note $F_{1}$ l'ensemble constitu\'e des deux intervalles ferm\'es $I_{1}^{(1)}$ et $I_{2}^{(1)}$. Le fait que les extr\'emit\'es de ces deux intervalles soient dans $P \backslash P_{0}$ permet de r\'eit\'erer le proc\'ed\'e sur chacun d'eux. \\
Ainsi \`a l'issu de la $\textrm{n}^{\textrm{i\`eme}}$ \'etape, on obtient un ensemble ferm\'e $F_{n}$ constitu\'e de $2^{n}$ intervalles ferm\'es de longueur inf\'erieure \`a $\dsp \frac{1}{3^{n}} (b_{0} - a_{0})$, en ayant retir\'e de nouveau $2^{n-1}$ intervalles ouverts $J_{k}^{(n)}$ ($1 \leq k \leq 2^{n-1}$) \`a extr\'emit\'es dans $P \backslash P_{0}$. On pose alors 
$$
F = \bigcap_{n=1}^{+\infty} F_{n}.
$$
$F$ est clairement un ensemble parfait contenu dans $P$. Il reste \`a v\'erifier que $F$ satisfait la condition de Carleson (C). On a $I_{0} \backslash F = \bigcup \limits_{n=1}^{+\infty} \bigcup \limits_{k=1}^{2^{n-1}} J_{k}^{(n)}$. Comme $F$ est de mesure nulle, il satisfait la condition de Carleson si et seulement si:
\begin{eqnarray}
\sum_{n=1}^{+\infty} \sum_{k=1}^{2^{n-1}} \big| J_{k}^{(n)} \big| \log{\frac{1}{\big| J_{k}^{(n)} \big|}} < +\infty. \label{somme}
\end{eqnarray}
Or chaque intervalle $J_{k}^{(n)}$ est de longueur inf\'erieure \`a $\dsp \frac{1}{3^{n-1}} (b_{0} - a_{0})$. Comme la fonction $\dsp x \mapsto x \log{\frac{1}{x}}$ est croissante sur $[ 0, e^{-1} ]$, pour $n$ assez grand, on a
\begin{eqnarray*}
\sum_{k=1}^{2^{n-1}} |J_{k}^{(n)}| \log{\frac{1}{|J_{k}^{(n)}|}} & \leq & 2^{n-1} \frac{b_{0} - a_{0}}{3^{n-1}} \log{\frac{3^{n-1}}{b_{0} - a_{0}}} \\
& \leq & (b_{0} - a_{0}) \big( (n-1) \log 3 - \log(b_{0} - a_{0}) \big) \Big( \frac{2}{3} \Big)^{n-1},
\end{eqnarray*}
et par cons\'equent (\ref{somme}) a bien lieu. 
\end{proof}

\begin{prop} \label{hypop}
Soit $E$ un ensemble ferm\'e non d\'enombrable du cercle unit\'e, alors la propri\'et\'e $P(0,t,\{ 1 \})$ n'est v\'erifi\'ee pour aucun r\'eel $t \geq 0$.
\end{prop}

\begin{proof}
Soit $E$ un ensemble ferm\'e non d\'enombrable du cercle unit\'e. D'apr\`es le lemme pr\'ec\'edent, $E$ contient un ensemble parfait $F$ qui v\'erifie la condition de Carleson. Soit alors $\mu$ une mesure continue \`a support inclu dans $F$ et $T_{\mu}$ l'op\'erateur d\'efini en (\ref{optimal1}). $T_{\mu}$ est une contraction dont le spectre est inclu dans $F$, et on montre en utilisant le lemme 2 de \cite{Atz1} que
\begin{eqnarray*} 
& & \big\| T_{\mu}^{-n} \big\| = O \big( e^{\e \sqrt{n}} \big) \, (n \rightarrow +\infty), \, \textrm{ pour tout } \e > 0. 
\end{eqnarray*}
Maintenant on conclut par des arguments bien connus (voir \cite{ERaZa}) que $\dsp \lim_{n \rightarrow +\infty} \big\| T_{\mu}^{-n} \big\| = +\infty$, et donc que $T_{\mu}$ n'est pas unitaire. Puis on d\'eduit du th\'eor\`eme 6.4 de \cite{Est1} que $T_{\mu}^{-n} = O(n^{t}) \, (n \rightarrow +\infty)$ n'est satisfait pour aucun r\'eel $t \geq 0$.  
\end{proof}

\ \\

\ \\
\ \\

\hspace*{-7mm} AGRAFEUIL Cyril \\
Cyril.Agrafeuil@math.u-bordeaux.fr \\
Universit\'e Bordeaux I \\
351, cours de la lib\'eration \\
33405 TALENCE C\'edex \\


\begin{thebibliography}{99}

\bibitem{ATaWi2} H. Alexander, B. A. Taylor and D. L. Williams, \emph{The interpolating sets for $\mathcal{A}^{\infty}$}, J. Math. Anal. Appl. \textbf{36} (1971), 556-566.

\bibitem{Atz1} A. Atzmon, \emph{Operators which are annihilated by analytic functions and invariant subspaces}, Acta Math. \textbf{144} (1980), 27-63.

\bibitem{BeGi} C. Bennett and J. E. Gilbert, \emph{Homogeneous algebras on the circle: I.-Ideals of analytic functions}, Ann. Inst. Fourier, \textbf{22} (1972), 1-19.

\bibitem{BourEVT4} N. Bourbaki, \emph{''El\'ements de math\'ematique. Espaces vectoriels topologiques ch. 1 \`a 5''}, Masson, 1981.

\bibitem{Carl} L. Carleson, \emph{Sets of uniqueness for functions regular in the unit circle}, Acta. Math. \textbf{87} (1952), 325-345.

\bibitem{Est1} J. Esterle, \emph{Uniqueness, strong form of uniqueness and negative powers of contractions}, Banach Center Publ. \textbf{30} (1994), 1-19.

\bibitem{ERaZa} J. Esterle, M. Rajoelina and M. Zarrabi, \emph{On contractions with spectrum contained in Cantor set}, Math. Proc. Camb. Phil. Soc. \textbf{117} (1995), 339-343.

\bibitem{EStZo1} J. Esterle, E. Strouse and F. Zouakia, \emph{Theorems of Katznelson-Tzafriri type for contractions}, J. Func. Anal. \textbf{94} (1990), 273-287.

\bibitem{EStZo2} J. Esterle, E. Strouse and F. Zouakia, \emph{Closed ideals of $A^{+}$ and the Cantor set}, J. reine angew. Math. \textbf{449} (1994), 65-79.

\bibitem{Hoff} K. Hoffman \emph{''Banach spaces of analytic functions''}, Prentic-Hall, Englewood Cliffs, 1962.

\bibitem{Kaha2} J. P. Kahane, \emph{Id\'eaux ferm\'es dans certaines alg\`ebres de fonctions analytiques}, Acte Table Ronde Inst. C.N.R.S. Montpellier, Lect. Notes Math. \textbf{336}, Springer Verlag (1973), 5-14.

\bibitem{Katz} Y. Katznelson, \emph{''An introduction to harmonic analysis''}, John Wiley ans Sons, 1968.

\bibitem{Kor} B. I. Korenblyum, \emph{Closed ideals in the ring $\mathcal{A}^{n}$}, J. Func. Anal. Appl. \textbf{6} (1972), 203-214.

\bibitem{Mat} A. L. Matheson, \emph{Closed ideals in rings of analytic functions satisfying a Lipschitz condition}, Lect. Notes Math. \textbf{604} Springer-Verlag (1976), 67-72.

\bibitem{NaFo} B. Sz. Nagy and C. Foias, \emph{''Analyse harmonique des op\'erateurs de l'espace de Hilbert''}, Akademiai Kiado, Budapest, 1967.

\bibitem{Rud2} W. Rudin, \emph{The closed ideals in an algebra of analytic functions}, Can. J. Math. \textbf{9} (1957), 426-434.

\bibitem{TaWi} B. A. Taylor and D. L. Williams, \emph{Ideals in rings of analytic functions with smooth boundary values}, Canad. J. Math. \textbf{22} (1970), 1266-1283.

\bibitem{TaWi2} B. A. Taylor and D. L. Williams, \emph{Zeros of Lipschitz functions analytic in the unit disc}, Mich. Math. J. \textbf{18} (1971), 129-139.

\bibitem{Zar1} M. Zarrabi, \emph{Synth\`ese spectrale dans certaines alg\`ebres de Beurling sur le cercle unit\'e}, Bull. Soc. Math. France \textbf{118} (1990), 241-249.

\bibitem{Zar2} M. Zarrabi, \emph{Contractions \`a spectre d\'enombrable et propri\'et\'e d'unicit\'e des ferm\'es d\'enombrables du cercle}, Ann. Inst. Fourier \textbf{43} (1993), 251-263.


\end{thebibliography}
\end{document}